\documentclass[12pt]{article}

\usepackage[T2A]{fontenc}
\usepackage[cp1251]{inputenc}
\usepackage[english,russian]{babel}
\usepackage[tbtags]{amsmath}
\usepackage{amsfonts,amssymb}

\newcommand{\dist}{\operatorname{dist}}
\newcommand{\spec}{\sigma}

\newcommand{\im}{\mathop{\rm Im}}
\newcommand{\re}{\mathop{\rm Re}}

\begin{document}

\title{\Large Preserving of the unconditional basis property under
non-self-adjoint perturbations of self-adjoint operators\thanks{This
work was supported by the Deutsche Forschungsgemeinschaft (DFG) and
by the Russian Foundation for Basic Research (RFBR).
}}

\author{\large A.\,K.~Motovilov, A.\,A.~Shkalikov}

\date{}

\maketitle

\begin{otherlanguage}{english}

\begin{abstract}
Let $T$ be a self-adjoint operator in a Hilbert space $H$  with
domain $\mathcal D(T)$. Assume that the spectrum of $T$ is confined
in the union of disjoint intervals $\Delta_k =[\alpha_{2k-1},\,
\alpha_{2k}]$, $k\in \mathbb{Z}$,  and
\begin{equation*}
\alpha_{2k+1}-\alpha_{2k} \geqslant b
|\alpha_{2k+1}+\alpha_{2k}|^p\quad \text{ for some }\, b>0,\,
p\in[0,1).
\end{equation*}%
Suppose that a linear operator $B$ in $H$ is $p$-subordinated to
$T$, i.e. $\mathcal D(B) \supset\mathcal D(T)$ and $\|Bx\| \leqslant
b'\,\|Tx\|^p\|x\|^{1-p}  +M\|x\| \text{\, for all } x\in \mathcal
D(T)$, with some $b'>0$ and  $M\geq 0$. Then the spectrum of the
perturbed operator $A=T+B$ lies in the union of a rectangle in
$\mathbb{C}$ and double parabola $P_{p,h} = \bigl\{\lambda \in
\mathbb{C}\, \bigl|\ \ |\mathop{\rm Im} \lambda| \leqslant
h|\mathop{\rm Re} \lambda|^p\bigr\}$, provided that $h>b'$. The
vertical strips $\Omega_k =\{\lambda\in\mathbb C\,\vert \ \,
|r_k-\re\,\lambda|\leqslant \delta r_k^p\}$, $\ r_k =
(\alpha_{2k}+\alpha_{2k+1})/2$,   belong to the resolvent set of
$T$, provided that $\delta <b -b'$ and $ |k|\geq N$ for $N$ large
enough. For  $|k|\geqslant N+1$,  denote by $\Pi_k$ the curvilinear
trapezoid formed by the lines $\re\,\lambda = r_{k-1}$,
$\re\,\lambda = r_{k}$, and the boundary of the parabola $P_{p,h}$.
Assume that $Q_0$ is the Riesz projection corresponding to the
(bounded) part of the spectrum of $T$ that lies outside
$\bigcup_{|k|\geqslant N+1}\Pi_k$. And let $Q_k$,  $|k|\geqslant
N+1$, be the Riesz projection for the part of the spectrum of $T$
confined within the curvilinear trapezoid $\Pi_k$. Main result of
the work consists in proving that the system of the invariant
subspaces $Q_k(H)$, $|k|\geqslant N+1$, together with the invariant
subspace $Q_0(H)$  forms an unconditional basis of subspaces in the
space $H$. We also prove a generalization of this theorem to the
case where any gap $(\alpha_{2k},\,\alpha_{2k+1})$,
$k\in\mathbb{Z}$, may contain a finite number of eigenvalues of $T$
with their total multiplicity bounded by a number $m\in\mathbb{N}$
independent of $k$.

\end{abstract}

\end{otherlanguage}
\vspace*{2cm}

\normalsize


\vskip 0.2cm


\vskip 0.2cm

\newpage

\begin{center}
\Large Сохранение безусловной базисности при несамосопряженных
возмущениях самосопряженных операторов$^*$
\end{center}
\vspace*{-2mm}

\begin{center}
\large А.\,К.~Мотовилов, А.\,А.~Шкаликов
\end{center}
\medskip

\bigskip

\begin{abstract}
Пусть $T$ --- самосопряженный оператор в гильбертовом пространстве
$H$ с областью определения $\mathcal D(T)$. Будем считать, что
спектр $T$ лежит в объединении непересекающихся интервалов $\Delta_k
=[\alpha_{2k-1},\, \alpha_{2k}]$, $k\in \mathbb{Z}$,  длины которых
подчиняются неравенствам
\begin{equation*}
\alpha_{2k+1}-\alpha_{2k} \geqslant b
|\alpha_{2k+1}+\alpha_{2k}|^p\quad \text{ при некоторых }\, b>0,\,
p\in[0,1).
\end{equation*}%
Предположим, что линейный оператор $B$, действующий в  $H$, является
$p$-под\-чи\-нен\-ным оператору $T$, т.е. $\mathcal D(B)
\supset\mathcal D(T)$ и $\|Bx\| \leqslant b'\,\|Tx\|^p\|x\|^{1-p}
+M\|x\|$ с некоторыми $b'>0$ и $M\geq 0$ при всех  $x\in \mathcal
D(T)$. Тогда при всяком $h>b'$ спектр возмущенного оператора $A=T+B$
лежит в объединении некоторого конечного прямоугольника в
$\mathbb{C}$ и двойной параболы $P_{p,h} = \bigl\{\lambda \in
\mathbb{C}\, \bigl|\ \ |\mathop{\rm Im} \lambda| \leqslant
h|\mathop{\rm Re} \lambda|^p\bigr\}$. Более того, если $\delta <b
-b'$, то можно указать  такое $N\in\mathbb{N}$, что при всех $|k|>N$
вертикальные полосы $\Omega_k =\{\lambda\in\mathbb C\,\vert \ \,
|r_k-\re\,\lambda|\leqslant \delta r_k^p\}$, $\ r_k =
(\alpha_{2k}+\alpha_{2k+1})/2$,  лежат в резольвентном множестве
$T$. При  $|k|\geqslant N+1$ обозначим через $\Pi_k$ криволинейные
трапеции, образованные прямыми $\re\,\lambda = r_{k-1}$,
$\re\,\lambda = r_{k}$ и границей параболы $P_{p,h}$. Пусть $Q_0$
--- проектор Рисса, отвечающий той (ограниченной) части спектра
оператора $T$, которая лежит вне $\bigcup_{|k|\geqslant N+1}\Pi_k$.
И пусть $Q_k$,  $|k|\geqslant N+1$, --- проектор Рисса для той части
спектра $T$, которая лежит в криволинейной трапеции $\Pi_k$.
Основным результатом работы является доказательство теоремы о том,
что система инвариантных подпространств $Q_k(H)$, $|k|\geqslant
N+1$, вместе с инвариантным подпространством $Q_0(H)$ образует
безусловный базис из подпространств в гильбертовом пространстве $H$.
Мы доказываем также обобщение этой теоремы на случай, когда в любой
из лакун $(\alpha_{2k},\,\alpha_{2k+1})$, $k\in\mathbb{Z}$, может
присутствовать конечный набор собственных значений $T$, суммарная
кратность которых не превышает некого числа $m\in\mathbb{N}$, не
зависящего от $k$.
\end{abstract}

\vfill

\noindent\rule{0.4\textwidth}{.1mm}

{\footnotesize\hspace*{-0.5cm}$^*$Работа выполнена при финансовой
поддержке Немецкого научно-исследовательского общества (DFG) и
Российского фонда фундаментальных исследований (РФФИ)}

\newpage

\begin{center}
$\S 1.$ {\bf Введение}
\end{center}

\medskip

Задача о спектральных свойствах операторов, которые представимы в
виде возмущений самосопряженных операторов, имеет давнюю историю.
Изложение разных аспектов этой теории для операторов с дискретным
спектром имеется в \cite{Sh1}.  В настоящей работе нас будет
интересовать вопрос о сохранении свойств базисности для подчиненных
возмущений  самосопряженных операторов  в общем случае без
предположения дискретности спектра невозмущенного оператора.
По-видимому, первым результатом  на эту тему была работа авторов
\cite{MSh}. Здесь мы получим существенное обобщение основного
результата  работы \cite{MSh}.

Одним из глубоких результатов, который можно рассматривать как итог
развития в 1951-1985 годах тематики, посвященной задаче о базисности
корневых векторов несамосопряженных операторов,   является теорема
Маркуса-Мацаева \cite[Гл. 1, Теорема 6.12]{Ma} (см. также
\cite{MaM1}). Приведем формулировку этой теоремы (для простоты
вместо нормальных операторов со спектром на конечном числе лучей мы
рассматриваем самосопряженные операторы).

{\bf  Теорема А. } {\sl  Пусть $T$ --- самосопряженный оператор с
дискретным спектром в гильбертовом пространстве $H$, а его
собственные значения $\{\mu_k\}$, занумерованные с учетом кратности,
подчинены оценке
\begin{equation}\label{1a}
|\mu_k| \geqslant C |k|^\beta, \qquad  C =\text{const}, \quad \beta
>0
\end{equation}
(здесь предполагается, что отрицательные собственные значения
нумеруются  в порядке их убывания отрицательными целыми числами, а
неотрицательные, в порядке их возрастания --- неотрицательными
индексами).  Пусть линейный (не обязательно замыкаемый) оператор $B$
задан на на области определения $\mathcal D(T)$  оператора $T$ и при
некоторых $b>0$, $0\leqslant p<1$ выполняется оценка
\begin{equation}\label{2}
\|Bx\| \leqslant b\|Tx\|^p\|x\|^{1-p},  \quad \ x\in \mathcal D(T).
\end{equation}
Тогда оператор $A=T+B$,  определенный на $\mathcal D(T)$, имеет
только дискретный спектр. При дополнительном условии $\beta^{-1}
\leqslant 1-p$ система его корневых векторов образует безусловный
базис со скобками в исходном гильбертовом пространстве $H$. }

Отметим, что условие \eqref{1a}  вместе с неравенством  $\beta^{-1}
\leqslant 1-p$ эквивалентно условию
\begin{equation}\label{0}
\overline{\lim}_{t\to\infty} \frac {n(t, T)}{t^{1-p}} <\infty,
\qquad\text {где }\ \,  n(t, T) = \sum_{|\mu_k|\leqslant t} 1\
\end{equation}
функция распределения собственных значений оператора $T$. Здесь же
отметим, что предложенное Маркусом и Мацаевым доказательство по
существу без изменений переносится на случай, когда верхний предел в
последнем соотношении меняется на нижний.

В работе \cite{Sh1}  было получено обобщение этой теоремы в
следующей форме.

{\bf  Теорема B. } {\sl  Пусть $T$ --- самосопряженный оператор с
дискретным спектром в гильбертовом пространстве $H$, а оператор $B$
является $p$-подчиненным оператору $T$,  то есть при некоторых
$p\in[0,1),$   выполняется оценка
 \begin{equation}\label{4}
\|Bx\| \leqslant b\|Tx\|^p\|x\|^{1-p}  +M\|x\|,  \quad \ x\in
\mathcal D(T), \quad b, M = \text{const}.
\end{equation}
Пусть $b'$ ---  точная нижняя грань, при которых выполняется оценка
\eqref{4}  с произвольными постоянными $M= M(b)$ (но при $p=0$
считаем, что постоянная $M$  в \eqref{4} равна нулю!)  и  при
некотором $b_1 > b'$ выполняется условие
 \begin{equation}\label{m}
\underline{\lim}_{\ r\to\infty}\,  n_\pm (r +b_1r^p) -n_\pm
(r-b_1r^p) =  m <\infty,
\end{equation}
где  $n_\pm (r)$ ---  функции распределения  собственных значений
оператора  $T$ на положительной и отрицательной осях соответственно.
Тогда  система корневых векторов оператора $A=T+B$  образует
безусловный базис со скобками. }

Отметим, что фигурирующее в Теореме B условие \eqref{m} существенно
слабее условия \eqref{1a} или  условия \eqref{0}.

Предположение о дискретности спектра оператора $T$  в Теореме B
существенно использовалось в ее доказательстве.  В частном, но
важном случае, когда порядок подчиненности $p=0$,   число $M$ в
условии \eqref{4} равно нулю, а вместо условия  \eqref{m} требуется
существование лакун $\Lambda_k = (\alpha_{2k},\alpha_{2k+1}) $ в
спектре $T$ c равномерно ограниченной снизу длиной, авторы
\cite{MSh} показали, что утверждение Теоремы B допускает обобщение
на операторы с непустым существенным и, в частности, с непрерывным
спектром. А именно, был получен следующий результат.

{\bf Теорема С.} {\sl Пусть $T$ --- самосопряженный оператор, спектр
которого лежит внутри  отрезков $\Delta_k = [\alpha_{2k-1},\,
\alpha_{2k}] \subset \mathbb R$, $k\in \mathbb Z$, причем
$$
  \inf_{k\in\mathbb Z} (\alpha_{2k+1} - \alpha_{2k}) = d >0
$$
(т.е. длины лакун в спектре между отрезками $\Delta_k$ всегда
$\geqslant d$). Пусть $B$ --- ограниченный (и, в общем,
несамосопряженный) оператор, причем $\|B\|=b^*<d/2$.  Тогда при
любом $\varepsilon \in (b^*,\, d/2)$ спектр возмущенного оператора
$A=T+B$  лежит в объединении  $\bigcup_k U_\varepsilon
\left(\Delta_k\right)$, где   $U_\varepsilon \left(\Delta_k\right)$
--- непересекающиеся между собой $\varepsilon$-окрестности отрезков
$\Delta_k$. Если $Q_k$ --- проекторы Рисса на инвариантные
подпространства оператора $A=T+B$, отвечающие частям спектра в
областях $U_\varepsilon \left(\Delta_k\right)$, т.е.
\begin{equation}\label{Riesz}
Q_k=-\frac{1}{2\pi i} \int_{\Gamma_k}
\left(A-\lambda\right)^{-1}\,d\lambda, \qquad \Gamma_k=\partial
U_\varepsilon \left(\Delta_k\right),
\end{equation}
то подпространства  $\mathcal{L}_k=Q_k(H)$, $k\in\mathbb Z$,
образуют безусловный  базис в исходном пространстве $H$. }

В этой работе мы получим сначала обобщение Теоремы С на случай
$p$-под\-чи\-нен\-ных возмущений при $p\in (0,1)$. А затем получим
еще более общий результат, который  является аналогом  Теоремы B
для операторов, которые могут  иметь  существенный спектр, но в
лакунах между отрезками существенного спектра имеется $\leqslant m$
собственных значений, где $m$ --- фиксированное число.

Перед формулировкой основных теорем введем следующие обозначения.
Обозначим через
\begin{equation}\label{par}
  P_{p,h} = \{\lambda \in \mathbb C\, \vert\ \ |\im \lambda| < h|\re \lambda|^p\}
\end{equation}
--- область, ограниченную двойной параболой. Часть этой области в правой
(левой) полуплоскости будем обозначать через $P^+_{p,h}\
(P^-_{p,h})$. Пусть $\Delta = [s,t]$ ---  отрезок на положительной
или отрицательной полуоси, $t>s$. Множество
\begin{equation}
\label{Updel} \mathcal U_{p\delta}(\Delta) = \{\lambda \in P_{p,
b'+\delta}\, \vert \ \ s  -\delta |s|^p  < \re\,\lambda  < t+\delta
|t|^p\}
\end{equation}
назовем $p\delta$-окрестностью отрезка $[s,t]$. Отметим, что в
определении $p\delta$-окрестности помимо чисел $p$  и $\delta$
участвует еще число $b'$, которое предполагается равным нижней грани
чисел $b$  в условии $p$-подчиненности \eqref{4}. Конечно, таким
образом определенная окрестность завит от этого числа $b'$, но для
удобства обозначений мы будем считать это число фиксированным и
указывать эту зависимость не будем. Всюду далее считаем, что в
случае полуограниченного оператора $T$ нумерация отрезков его
спектра проводится индексами $k\geqslant 1$, а в общем случае ---
целыми индексами  $k\in \mathbb Z$,  причем для этого случая
принимается $\alpha_0 \leqslant 0, \alpha_1 >0$.

%
%

Сформулируем основные результаты работы.
\smallskip

{\bf Tеорема 1.} {\sl Пусть $T$ --- самосопряженный оператор в
гильбертовом пространстве $H$  с областью определения $\mathcal
D(T)$,  а линейный оператор $B$  определен на  области $\mathcal
D(T)$ и является $p$-подчиненным оператору $T$  (в смысле
неравенства \eqref{4})  с нижней гранью $p$-подчиненности $b'$ (как
и ранее  при $p=0$  в неравенстве \eqref{4} полагаем $M=0$).
 Предположим, что  спектр  оператора $T$ заключен внутри
непересекающихся отрезков $\Delta_k =[\alpha_{2k-1},\,
\alpha_{2k}]$, а длины лакун $\Lambda_k = (\alpha_{2k},\,
\alpha_{2k+1})$ между отрезками $\Delta_k$ таковы,
что\footnote{Условие \eqref{gaps} в  Теореме 1  эквивалентно тому,
что лакуны $\Lambda_k$ содержат интервалы $(r_k -b_1 |r_k|^p,\, r_k
+b_1 |r_k|^p)$, где $r_k =(\alpha_{2k}+\alpha_{2k+1})/2$ ---
середины лакун $\Lambda_k$.}
\begin{equation}\label{gaps}
\alpha_{2k+1}-\alpha_{2k} \geqslant 2^{1-p}\,\, b_1
|\alpha_{2k+1}+\alpha_{2k}|^p\quad \text{при некоторых}\ \, b_1>0,\,
p\in[0,1),
\end{equation}
причем $b_1 > b'$. Тогда при любом $\delta \in (0, b_1 -b')$
найдется число $N\in \mathbb N$, такое, что  спектр   оператора
$A=T+B$    лежит в объединении некоторого прямоугольника $\mathcal
U_0$ и взаимно непересекающихся $p\delta$-окрестностей $ \mathcal
U_k =\mathcal U_{p\delta}(\Delta_k)$ отрезков
 $\Delta_k$  при $|k| \geqslant N$.  Если
\begin{equation*}
 Q_k=-\frac{1}{2\pi i} \int_{\partial\mathcal U_k}
\left(A-\lambda\right)^{-1}\,d\lambda,\qquad k=0, \pm N, \pm (N+1),
\dots,
\end{equation*}
--- проекторы Рисса,  отвечающие
частям спектра оператора $A$, расположенным внутри $\mathcal U_k$,
то инвариантные подпространства $\{Q_k(H)\}$, $|k|\geqslant N$,
вместе с $Q_0(H)$ образуют безусловный базис из подпространств в
гильбертовом пространстве $H$. }
\smallskip

{\bf Tеорема 2.} {\sl Пусть выполнены условия Теоремы 1, но в
ослабленном виде: допускается, что вне отрезков $\Delta_k$,
содержащих существенный спектр оператора $T$, имеется (только)
дискретный спектр, причем в каждой лакуне $\Lambda_k$  содержится
$\leqslant m$  собственных значений оператора $T$ с учетом кратности
и число $m$ не зависит от $k$. Тогда спектр оператора $A=T+B$
асимптотически лежит внутри двойной параболы $  P_{p,h}$, $h>b'$, и
при любом  $\delta\in (0, b_1-b')$ вертикальные полосы
$$
\Omega_k =\{\lambda\in \mathbb C\, \vert\ \,  r_k -\delta r_k^p
\leqslant \re\, \lambda \leqslant r_k+\delta r^p_k\}, \qquad r_k =
(\alpha_{2k}+ \alpha_{2k+1})/2,
$$
не могут содержать никаких иных точек спектра оператора $A$ кроме
собственных значений c суммарной алгебраической кратностью
$\leqslant Cm$, где $C$ не зависит от $k$. Найдутся такие  числа
$r'_k\in \Omega_k \bigcap \mathbb R$, что границы криволинейных
трапеций $\Pi_k$, образованных прямыми $\re \lambda = r'_k$,  $\re
\lambda = r'_{k+1}$ и параболой $P_{p,h}$, не содержат точек спектра
оператора $A=T+B$, а проекторы Рисса
\begin{equation*}
 Q_k=-\frac{1}{2\pi i} \int_{\partial\Pi_k}
\left(A-\lambda\right)^{-1}\,d\lambda,\qquad k=\pm N, \pm (N+1),
\dots,
\end{equation*}
таковы, что  подпространства $Q_k(H)$, $|k|\geqslant N$, вместе с
$Q_0(H)$ образуют безусловный базис из подпространств в гильбертовом
пространстве $H$. Здесь $Q_0$ ---  проектор Рисса, отвечающий той
части спектра оператора $A$, которая не входит в криволинейные
трапеции $\Pi_k$, $|k| \geqslant N$. }


\begin{center}
$\S 2.$ {\bf  Определения и вспомогательные предложения.}
\end{center}

Напомним известные понятия базиса и безусловного базиса (см.
например, \cite [гл. 6]{GK}. Последовательность векторов
$\{\psi_k\}_{k=1}^\infty$ пространства $H$ называется {\it базисом}
этого пространства, если каждый вектор $x\in H$ разлагается
единственным образом в сходящийся по норме $H$ ряд
\begin{equation}\label{ser}
x=\sum_{k=1}^\infty  c_k\psi_k,
\end{equation}
где $c_k$ --- числовые коэффициенты. Базис называется {\it
безусловным}, если он остается базисом после любой перестановки
векторов последовательности. Система $\{\psi_k\}_{k=1}^\infty$
называется {\it базисом со скобками}, если для любого $x\in H$ ряд
\eqref{ser} сходится по норме $\mathcal H$ после некоторой
расстановки скобок, не зависящей от $x$. Если система остается
базисом после любой перестановки  наборов ее векторов, которые
отвечают членам ряда, заключаемым в скобки, то такая система
называется {\it безусловным базисом со скобками} или {\it базисом
Рисса со скобками}.

Далее будет идти речь о базисности и безусловной базисности
последовательности подпространств, о чем удобно говорить на языке
проекторов (см. \cite[$\S 6$]{Sh1}). Используются следующие
определения:
 \begin{enumerate}
 \item [({\it a})]
Система ограниченных проекторов $ \{Q_j\}_1^\infty$ полна, если
равенства
$$
(Q_jx, y) =0  \quad \forall x\in H, \ \ j=1,2,\dots,\ \
\text{влекут}\ \, y=0.
$$
\item [({\it b})]
Система ограниченных проекторов $ \{Q_j\}_1^\infty$  минимальна,
если $ Q_jQ_k =\delta_{kj} Q_j$.
\item [({\it c})]
Пусть $ \{Q_j\}_1^\infty$ --- минимальная система ограниченых
проекторов в $H$. Система подпространств $\mathcal L_j = Q_j(H),
j=1,2,\dots,$ называется базисом (безусловным базисом)  в $H$,  если
ряд $\sum_j Q_j$  сходится (сходится после любой перестановки
индексов в сумме) в сильной операторной топологии к единичному
оператору.
\end{enumerate}

Далее будем использовать следующий результат \cite[$\S 6$]{Sh1}.
\smallskip

{\bf Предложение 1.} {\sl Пусть $\{Q_j\}_{j=1}^\infty$  --- система
ограниченных проекторов в $H$. Система подпространств $\mathcal L_j
=Q_j(H), \ j=1,2, \dots,$  является безусловным базисом из
подпространств в $H$ тогда и только тогда, когда система проекторов
$\{Q_j\}_{j=1}^\infty$ полна, минимальна и
\begin{equation}\label{converge}
\sum_{j=1}^\infty |(Q_jx,x)| < \infty \qquad \forall\, x\in H.
\end{equation}
}

Приведем утверждение о локализации спектра оператора оператора $A=
T+B$   при $p$-подчиненных возмущениях.
\smallskip

{\bf Предложение 2.} {\sl Пусть $T$ --- самосопряженный
положительный (или полуограниченный) оператор, а возмущение $B$
является $p$-подчиненным оператору $T$ с точной верхней гранью $b'$.
Тогда спектр  оператора  $A=T+B$ асимптотически локализован в
параболической области (при $p=0$ полуполосе)
\begin{equation}\label{Pi}
P^+_{p,\, h} =\{\lambda\in\mathbb{C}\, | \ \, \re \lambda \geqslant
0,\ \,  |\im\lambda| \leqslant h(\re \lambda)^p\}, \quad h>b'.
\end{equation}
При этом в правой полупоскости вне параболы $P^+_{p,\, h}$ при любом
$\varepsilon >0$  асимптотически выполняются оценки (при любом
$\varepsilon\in (0, h-b')$)
\begin{equation}\label{resS}
 \|B(T-\lambda)^{-1}\| \leqslant \frac{b'+\varepsilon}{h} \qquad |\lambda| \geqslant R = R(\varepsilon),
 \end{equation}
и
 \begin{equation}\label{resA}
\|(A-\lambda)^{-1}\| \leqslant
\left(\frac{b'+\varepsilon}{h-b'}\right)\, \frac 1{|\im \lambda|},
\qquad |\lambda| \geqslant R,
 \end{equation}
а в левой полуплоскости вне круга радиуса $\gg 1$ выполняется оценка
оценка \eqref{resS} и оценка  $\|(A-\lambda)^{-1}\| \leqslant
C|\lambda|^{-1}$.  Если $T$ самосопряжен, но не полуограничен, то
обе оценки \eqref{resS} и \eqref{resA} выполняются вне области,
ограниченной двойной параболой $P_{p,\, h}$.%
}%
\smallskip

При выполнении условия \eqref{2} это утверждение доказано в  работе
\cite{MaM2}. При выполнении более общего условия  $p$-подчиненности
\eqref{4} доказательство Предложения~1 имеется в \cite[$\S 3$]{Sh1}.
\smallskip

{\bf Предложение 3.} {\sl Пусть $\{s_j\}_1^\infty$ и
$\{\lambda_j\}_1^\infty$ --- $s$-числа и собственные значения
ядерного оператора $K$, занумерованные с учетом алгебраической
кратности. Тогда
\begin{equation}\label{K}
(1-\mu K)^{-1} = \frac{F(\mu)}{D(\mu)}, \quad \text{где}\ \, D(\mu)
= \prod_{j=1}^\infty \left(1-\mu\,\lambda_j\right).
\end{equation}
Произведение $D(\mu)$  сходится и является целой функцией переменной
$\mu\in\mathbb{C}$ порядка $\leqslant 1$ и минимального типа при
порядке 1. Оператор-функция $F(\mu)$ также целая и допускает оценку
$$
\|F(\mu)\| \leqslant \prod_{j=1}^\infty \left( 1+ |\mu| s_j\right).
$$
Тем самым $F(\mu)$  также порядка $\leqslant 1$  и минимального типа
при порядке 1.}%
\smallskip

{\it Доказательство}  этого важного утверждения имеется в \cite[гл.
3]{GK}.

Нам понадобится также одно  предложение из теории  функций.
\smallskip

{\bf Предложение 4.}
{\sl Пусть $f(\lambda)$ голоморфная ограниченная функция в $\lambda$
прямоугольнике
\begin{equation}\label{Pi}
\Pi\ = \{\lambda: \ \, |\re\lambda| <a, \ \, |\im\lambda| <c\}.
\end{equation}
Фиксируем $\tau \in (0,1)$, положим $a' =a(1-\tau),\ \,
c'=c(1-\tau)$ и  обозначим через $\Pi'$ прямоугольник,  который
определен равенством \eqref{Pi}  с заменой $a$ и $c$ на $a'$ и $c'$,
соответственно. Обозначим
$$
 {M}= \sup_{\lambda\in \Pi} f(\lambda),\quad {M}' = \sup_{\lambda\in \Pi'} f(\lambda).
$$
Тогда для числа нулей $n_f(\Pi')$  функции $f(\lambda)$ в
прямоугольнике $\Pi'$  справедлива оценка
\begin{equation}\label{n_f}
n_f(\Pi') \leqslant C(\ln M-\ln M').
\end{equation}
Для любого $\varepsilon\in (0,1)$  найдутся круги с общей суммой
радиусов $\leqslant \varepsilon$, такие, что  вне объединения
$\mathcal E$ этих кругов выполняется оценка снизу
\begin{equation}\label{below1}
\ln|f(\lambda)| \geqslant С \ln\varepsilon\ (\ln M -\ln M') + \ln
M', \quad \lambda\in\Pi' \setminus \mathcal E.
\end{equation}
Постоянная $C$ в обеих оценках зависят только от $\tau$ и отношения
$a/c$ и не зависит от $f$. При  одновременном увеличении или
уменьшении сторон прямоугольников $\Pi$ и $\Pi'$ в $\eta$ раз оценка
\eqref{below1}  сохраняется с той же постоянной $C$, но
исключительное множество кругов нужно брать  с общей суммой радиусов
$\eta\varepsilon$. }

Варианты   сформулированных в этом предложении утверждений в другой
форме можно найти, например, в книге \cite{Lev}. В сформулированном
виде первое  утверждение  имеется в работе \cite[Лемма 1.1 ]{MaM2},
а первое и второе  утверждение в работах
 \cite{Sh4} и \cite[Лемма 2.1]{Sh6}.

  \begin{center}
$\S 3.$ {\bf   Доказательство основных теорем.}
\end{center}

Поскольку Теорема 2 является обобщением Теоремы 1, мы не будем
проводить их доказательства отдельно. Сначала проведем основную
часть доказательства в условиях Теоремы 1, затем укажем на
изменения, которые нужно сделать в общем случае, а на заключительном
этапе докажем полноту проекторов Рисса  в условиях Теоремы 2.
Доказательство разобьем  на несколько этапов. На некоторых этапах
доказательство проводится почти так же, как  в случае Теореме 6.6 из
работы \cite{Sh1}. Однако для удобства читателя мы приводим
доказательство полностью. Всюду далее через $\spec(Z)$ обозначается
спектр $Z$.

{\it Шаг 1.} Упростим максимально  задачу  и будем считать, что
оператор $T=T^*$ положителен,    а  лакуны
$\Lambda_k:=(\alpha_{2k},\alpha_{2k+1})$  в спектре невозмущенного
оператора $T$ не содержат  собственных значений.
Позже мы избавимся от этих предположений.  Доказательство проведем
на основе следующих представлений, справедливых при $\lambda \notin
\sigma(T) \cup \sigma(A)$:
\begin{equation}\label{65}
(A-\lambda)^{-1} = (T-\lambda)^{-1} - G(\lambda), \qquad G(\lambda)
= (A-\lambda)^{-1}B(T-\lambda)^{-1}
\end{equation}
и
\begin{equation}\label{66}
G(\lambda)  =(T-\lambda)^{-1} M(\lambda)B(T-\lambda)^{-1}, \quad
M(\lambda) = (1+S(\lambda))^{-1}, \ \ S(\lambda) =
B(T-\lambda)^{-1}.
\end{equation}
Согласно Предложению 2, спектр оператора $A= T+B$ лежит в
объединении  круга достаточно большого радиуса и  параболической
области $P_{p, h}$,  если $h>b'$, а вне этого объединения выполнена
оценка \eqref{resA}. Напоминаем, что $b'$ есть точная нижняя грань
постоянных $b$,  при которых выполняется оценка $p$-подчиненности
\eqref{4}.  Далее считаем $h > 2b'$,  тогда в силу Предложения 2
существует такое $r_0$, что
$$
\|S(\lambda)\|\leqslant 1/2, \ \ \  \|M(\lambda)\| \leqslant 2,
\qquad\text{когда \,\,} \lambda\notin P_{p, h}, \ \, |\lambda| >r_0.
$$

Ниже на шаге 4 мы покажем, что если интервалы $\Lambda_j=
(\alpha_{2j},\, \alpha_{2j+1})$  не содержат собственных значений
оператора $T$, то  при $|j|\geqslant N$ (где $N$ --- достаточно
большое число) на вертикальных прямых $\re\, \lambda =r_j=
(\alpha_{2j+1}+ \alpha_{2j})/2 $  справедлива оценка
$\|S(\lambda)\|\leqslant 1-\varepsilon$  при некотором $\varepsilon
>0$,  а потому оператор-функция $M(\lambda)$  равномерно ограничена
на этих прямых. В частности, эта оператор-функция равномерно
ограничена на вертикальных отрезках $\gamma_j$,  которые проходят
через точки $r_j$ и соединяют кривые
$$
     \Gamma_\pm = \{\lambda \in \mathbb C \, \vert \ \im\, \lambda =
\pm h|\re\,\lambda|^p\}
$$
--- границы области $P_{p, h}$  в верхней и нижней полуплоскостях.

Обозначим через  $Q_j$ и $P_j$ проекторы Рисса на части спектра
операторов $A$ и $T$, которые заключены внутри криволинейных
трапеций $\Pi_j \subset P_{p, h}$, ограниченных вертикальными
отрезками $\gamma_{j-1}$ и $\gamma_j$. Сумма ортогональных
проекторов $P_j$  безусловно сходится к единичному оператору.
Система проекторов $\{Q_j\}$, отвечающих за непересекающиеся части
спектра, минимальна в  $H$ для любого оператора. Поэтому в силу
представления \eqref{65} для доказательства теоремы достаточно
показать, что система подпространств $Q_j(H)$, инвариантных
относительно оператора $A$, полна в $H$  и  при всех $x\in H$
сходится  ряд
\begin{equation}\label{67}
  \sum_{j=N}^\infty |(Q_j x, x)|, \qquad Q_j = -\frac 1{2\pi i}
\int_{\Gamma_j} (A-\lambda)^{-1}\, d\lambda,
\end{equation}
где $\Gamma_j$ --- ориентированные границы криволинейных трапеций
$\Pi_j$.

Задача о полноте инвариантных подпространств оператора $A$  будет
положительно решена  на заключительном шаге 6. Если равномерная
ограниченность опе\-ра\-тор-функции $M(\lambda)$  на вертикальных
отрезках $\gamma_j$ уже доказана,  то из представления \eqref{66}
получаем, что для завершения доказательства теоремы достаточно
доказать сходимость интегралов
\begin{equation}\label{Gam}
\int_{\Gamma_\pm} \|(T-\lambda)^{-1}x\|\ \| S(\lambda)x\|\
|d\lambda| \ < \infty, \qquad x\in H,
\end{equation}
и ряда
\begin{equation}\label{Ser1}
\sum_{j=1}^\infty \int_{\gamma_j} \|(T-\lambda)^{-1}x\|\ \|
S(\lambda)x\|\ |d\lambda| \ < \infty.
\end{equation}
Итак, наш план состоит в следующем.  На шаге 2 мы оценим
оператор-функцию $M(\lambda)$ на вертикальных прямых $\re\,\lambda =
r_j$  при условии, что лакуны $\Lambda_j = (\alpha_{2j},
\alpha_{2j+1})$  не содержат собственных значений оператора $T$.
Более того, равномерная оценка для $M(\lambda)$ будет сохраняться в
вертикальных полосах
\begin{equation}\label{Om}
\Omega_j = \{\lambda \in \mathbb C\, \vert \ \,  r_j-\delta r_j^p <
\re\,\lambda < r_j+\delta r_j^p\}
\end{equation}
при любом  фиксированном $\delta\in(0, b_1- b')$, где число $b_1$
берется из условия \eqref{gaps}. Далее мы докажем сходимость
интегралов \eqref{Gam} и сходимость ряда \eqref{Ser1}, где
$\gamma_j$ произвольные вертикальные отрезки в полосах $\Omega_j$,
соединяющие кривые $\Gamma_{\pm}$. На шаге 5 с помощью приема
искусственной лакуны мы покажем, что при выполнении условия
\eqref{m} с произвольным числом $m\in \mathbb N$ в полосе $\Omega_j$
можно выбрать вертикальные отрезки  $\gamma_j'$, на которых нужные
оценки сохраняются. Тем самым, доказательство будет завершено для
случая $T=T^* >0$.  Изменения, которые нужно провести для
доказательства в случае, когда оператор $T$  не является
полуограниченным, будут отражены на шаге~7.

{\it Шаг 2.} Оценим оператор-функцию $S(\lambda)$   на прямой
$\re\,\lambda =r_j$ при условии, что интервал $\Theta_j: = (r_j- b_1
r_j^p,\,  r_j + b_1 r_j^p)$  не содержит собственных значений
оператора $T$. Далее всюду, где это удобно, опускаем индекс $j$ и
вместо $r_j$ пишем $r$.

Пусть  $\lambda = \eta +i\tau$, \, $\eta,\tau\in\mathbb{R}$.
Согласно спектральной теореме для самосопряженного оператора $T$
имеем
   \begin{multline}\label{TT}
 \|T(T-\lambda)^{-1}\| = \sup_{t\in \sigma(T)} \left|\frac t{t-\lambda}\right| =
 \left( \inf_{t\in \sigma(T)} \left| 1-\frac\lambda{t}\right|\right)^{-1}
\leqslant  \\
  \left( \inf_{t\in \sigma(T)} \left|
1-\frac\eta{t}\right|\right)^{-1}, \quad \lambda\not\in\spec(T).
\end{multline}
Так как $\sigma(T)\cap \Theta_j = \emptyset$,  а функция $|1-\eta
t^{-1}|$ убывает при $t\in (0, \eta)$ и   возрастает при $t\in
(\eta, \infty)$,  то правая часть \eqref{TT}  на прямой
$\re\,\lambda = r =r_j$ не превосходит  максимума значений этой
функции на концах интервала $\Theta_j$,  то есть не превосходит
величины
$$
\max \left(\left| 1- \frac{r}{r-b_1r^p}\right|^{-1}, \ \left( 1
-\frac{r}{r+b_1r^p}\right)^{-1}\right)= \frac{r+b_1r^p}{b_1r^p}.
$$
В силу условия \eqref{gaps}  середина лакуны $\Lambda_j$  удалена от
концов отрезков, содержащих спектр оператора $T$, на расстояние
$\geqslant b_1 r^p$. Тогда, используя очевидную оценку
$\|(T-\lambda)^{-1}\|\leqslant (b_1 r^p)^{-1}$  на прямой
$\re\,\lambda =r= r_j$ и оценку \eqref{4} получаем (далее $o(1)$
есть бесконечно малая величина при $r =r_j\to\infty$)
 \begin{multline*}
\|S(\lambda)\| \leqslant b \left(\frac{r+ b_1r^p}{b_1r^p}\right)^p
\left( \frac 1 {b_1r^p}\right)^{1-p} + \frac M{b_1r^p} = b
\left(\frac{r+ b_1r^p}{b_1r^p}\right)^p \left( \frac 1
{b_1r^p}\right)^{1-p} \left(1+o(1)\right)\\
=\frac b{b_1}\left(1+o(1)\right) <1,
\end{multline*}
если $b\in (b', b_1)$.  Очевидно, последняя оценка сохраняется для
$\lambda$ в полосе $\Omega_r = \Omega_j(\delta)$, определяемой
\eqref{Om},  если $\delta\in (0, b_1-b)$ и $r_j$ достаточно велико.
Следовательно, спектр возмущенного оператора $A= T+B$ асимптотически
лежит в $p\delta$-окрестностях $\mathcal U_{p\delta}(\Delta_j)$ (см.
определение \eqref{Updel} в $\S$ 1).  В частности, при $\delta =
(b_1 -b')/2$ и $b= b'+ \delta/4$ получаем
\begin{equation}\label{Sl}
 \|S(\lambda)\| \leqslant \frac b{(b_1-\delta)} = \frac {2b}{b_1+
b'}, \quad \| M(\lambda)\|\leqslant \frac{b_1+ b'}{b_1+b' -2b}=
\frac{2(b_1+ b')}{(b_1-b')}, \qquad \lambda\in \Omega_j, \ \lambda
\gg 1.
 \end{equation}

{\it Шаг 3.} Оценим интеграл \eqref{Gam} по кривой $\Gamma_+$.
Оценка по $\Gamma_-$  аналогична. Из условия $p$-подчиненности
\eqref{4}   получаем
\begin{equation}\label{Slam}
\| S(\lambda)x\| \leqslant b\| T(T-\lambda)^{-1}x\|^p\
\|(T-\lambda)^{-1}x\|^{1-p} + M\|(T-\lambda)^{-1}x\|.
\end{equation}
При $p=0$ ситуация сильно упрощается,  так как второе слагаемое
равно нулю, а в первом слагаемом нужно работать только со вторым
множителем (соответствующие оценки при $p=0$  имеются в работе
\cite{MSh}). При $p>0$  интеграл в \eqref{Gam}  мажорируется
величиной
\begin{equation*}
b\,I_1+M\,I_2,
\end{equation*}
где
\begin{align}
I_1 &=\int_{\Gamma_+}\,  \left(\|T(T-\lambda)^{-1}x\|^p\
|\lambda|^{\frac{(p-2)p}2}\right)\
\left( |\lambda|^{\frac{(2-p)p}2}\|(T-\lambda)^{-1}x\|^{2-p}\right)\ |d\lambda|,\\
I_2 &=\int_{\Gamma_+}\, \|(T-\lambda)^{-1}x\|^{2}\, |d\lambda|.
\end{align}
Легко видеть, что при $\lambda\to\infty$, $\lambda\in\Gamma_+$,
подынтегральная функция в  $I_2$ являетcя бесконечно малой более
высокого порядка нежели подынтегральная функция в $I_1$. Поэтому
достаточно убедиться в сходимости интеграла $I_1$. Мы оценим этот
интеграл по неравенству Гёльдера (в качестве гёльдеровских
сопряженных чисел берем $q= 2/p$ и  $q' = 2/(2-p)$:
\begin{equation}\label{Gel1}
I_1\leqslant \left( \int_{\Gamma_+} |\lambda|^{p-2}\
\|T(T-\lambda)^{-1}x\|^2\ |d\lambda|\right)^{p/2} \left(
\int_{\Gamma_+} |\lambda|^p\ \|(T-\lambda)^{-1}x\|^2\, |d\lambda|
\right)^{(2-p)/2}.
\end{equation}
\tolerance=1000 Обозначим   через $E(t)$  спектральную функцию
оператора $T$  и положим $e(t) =(E(t)x,x).$ Пусть $\xi\geqslant 1$ и
$\lambda =\xi +ih\xi^p\,\,(\in \Gamma_+)$. Очевидно, при
$\xi\to+\infty$ имеем
$$
|\lambda|^p = \xi^p \left(1+ o(1)\right), \ \, |d\lambda| =
d\xi\left(1+o(1)\right),\ \, |t-\lambda|^2 = |t-\xi|^2
+ih^2\xi^{2p}.
$$
Поэтому сходимость интегралов в \eqref{Gel1} эквивалентна сходимости
следующих интегралов
\begin{equation*}
\int_1^\infty\hspace*{-0.5em} d \xi \int_0^\infty
\hspace*{-0.5em}d\, e(t)\frac{\xi^p\,}{(t-\xi)^2 + h^2\xi^{2p}}
=\int_1^\infty F_1(t)\, d\, e(t),\, \text{ где }\, F_1(t)
=\int_1^\infty  \frac{\xi^p}{(t-\xi)^2 + h^2\xi^{2p}} \, d \xi,
\end{equation*}
\medskip
\begin{equation*}
\int_1^\infty\hspace*{-0.5em} d\xi \int_0^\infty\hspace*{-0.5em} d\,
e(t) \frac{t^2\, \xi^{p-2}}{(t-\xi)^2 + h^2\xi^{2p}}  =\int_1^\infty
F_2(t)\, d\, e(t),\, \text{ где }\, F_2(t) =\int_1^\infty
\frac{t^2\,\xi^{p-2}}{(t-\xi)^2 + h^2\xi^{2p}} \, d \xi.
\end{equation*}
Поскольку функция $e(t)$ монотонна и $0\leqslant e(t) \leqslant
\|x\|^2$, выписанные интегралы сходятся, если  функции $F_1$  и
$F_2$  ограничены. Заметим, что
$$
F_1(t) \geqslant F_2(t), \ \text{если}\,  t\geqslant \xi, \
\text{и}\ F_1(t) \leqslant F_2(t), \ \text{если}\, 0\leqslant
t\leqslant \xi.
$$
Отсюда заключаем, что, в частности, при $t\geqslant 1$
$$
F_1(t) + F_2(t) \leqslant 2 \int_1^t \frac{t^2\,
\xi^{p-2}}{(t-\xi)^2 + h^2\xi^{2p}}\, d\xi +2 \int_t^\infty
\frac{\xi^{p}}{(t-\xi)^2 + h^2\xi^{2p}}\, d\xi =:2F_3(t) +2F_4(t).
$$
Оценка последних интегралов очевидна. Например, функция $F_3(t)$ при
$t^p<t/2$ (т.е. при $t>2^{1/(1-p)}$) оценивается следующим образом:
\begin{multline*}
F_3(t) \leqslant \left(\int_1^{t/2} + \int_{t/2}^{t-t^p} +
\int_{t-t^p}^t\right)
\frac{t^2 \, \xi^{p-2}}{(t-\xi)^2 +h^2\xi^{2p}}\, d\xi \\
\leqslant \int_1^{t/2}\, \frac{t^2\,\xi^{p-2}}{(t-\xi)^2}\, d\xi +
\int_{t/2}^{t-t^p} \,\frac{t^p}{(t-\xi)^2}\, d\xi+\frac{1}{h^2}
\int_{t-t^p}^t\, \frac{t^p}{\xi^{2p}}\,d\xi < \frac 4{1-p}
+1+\frac{4^p}{h^2}.
\end{multline*}
Аналогично оценивается интеграл $F_4(t)$. Тем самым оценка
\eqref{Gam} доказана.

{\it Шаг 4.} Докажем сходимость ряда \eqref{Ser1}.  Так как
интервалы $\Delta_n$  и $\Delta_{n+1}$ не пересекаются, то
$r_{n+1}\geqslant r_n +b_1(r_n^p+ r^p_{n+1})$. Поэтому при некотором
$n_0\in\mathbb{N}$
$$
r_{n+1}^{1-p} \geqslant r_n^{1-p}\left(1+ 2b_1\,
r_n^{p-1}\right)^{1-p} > r_n^{1-p}+ b_1(1-p), \quad \text{если}\
n\geqslant n_0.
$$
Следовательно, $r_{n+1}^{1-p} - r_n^{1-p}> b_1(1-p)$, если
$n\geqslant n_0$. Тогда по индукции при $j\geqslant n$ получаем
\begin{equation}\label{1-p}
r_j^{1-p} - r_n^{1-p}> c(j-n), \qquad c= b_1(1-p), \ \, j\geqslant
n\geqslant n_0.
\end{equation}
Поэтому
\begin{equation*}
r_j-r_n \geqslant r_j^p (r_j^{1-p} -r_n^{1-p}) \geqslant c\, r_j^p\,
(j-n), \quad j\geqslant n\geqslant n_0.
\end{equation*}
Меняя $j$ и $n$  местами, получаем
\begin{equation}\label{r}
|r_j-r_n| \geqslant c\max (r_j^p, \, r_n^p)\, |j-n|.
\end{equation}

Подынтегральная функция в \eqref{Ser1} с учетом неравенства
\eqref{4}  оценивается величиной
$$
\|T(T-\lambda)^{-1} x\|^p\, \|(T-\lambda)^{-1} x\|^{2-p} \,  + \,
M\|(T-\lambda)^{-1} x\|^2.
$$
Здесь главную роль играет первое слагаемое, поскольку второе
слагаемое  при $\lambda\to\infty$  является бесконечно малой более
высокого порядка. Поэтому  интегралы в  \eqref{Ser1} оцениваются
сверху через интегралы
\begin{equation}
\label{Intj} \int_{\gamma_j}\, \left(\|T(T-\lambda)^{-1} x\|^p\,
r_j^{\frac{p(p-2)}2}\right)\,\left(  r_j^{\frac{p(2-p)}2}  \,
\|(T-\lambda)^{-1} x\|^{2-p}\right)\, |d\lambda|, \quad \re\,\lambda
= r_j.
\end{equation}
Теперь для  оценки  ряда \eqref{Ser1}  применим  сначала неравенство
Гёльдера для интегралов, а затем неравенство Гёльдера для сумм. В
результате  для ряда \eqref{Ser1} придем к оценке сверху величиной
\begin{align}
\nonumber &\qquad\sum_{j=1}^\infty\ \left(\int_{\gamma_j}\,
r_j^{p-2}\| T(T-\lambda)^{-1}x\|^2\, |d\lambda|\right)^{p/2}
\left(\int_{\gamma_j}\, r_j^p\|(T-\lambda)^{-1}x\|^2\, |d\lambda|
\right)^{(2-p)/2}\\
\label{Prod} &\quad\leqslant \left(\sum_{j=1}^\infty
\int_{\gamma_j}\, r_j^p\|(T-\lambda)^{-1}x\|^2\,
|d\lambda|\right)^{(2-p)/2} \left(\sum_{j=1}^\infty
\int_{\gamma_j}\, r_j^{p-2}\| T(T-\lambda)^{-1}x\|^2\,
|d\lambda|\right)^{p/2}.
\end{align}
Как и прежде, будем использовать обозначение $e(t) = (E(t)x,\, x)$,
где $E(t)$ --- спектральная функция оператора $T$. При $\lambda \in
\gamma_j$ имеем
$$
\|(T-\lambda)^{-1}x\|^2 = \int_0^\infty \frac{de(t)}{(t- r_j)^2 +
\tau^2} \leqslant \int_0^\infty \frac{de(t)}{(t- r_j)^2}, \quad
\lambda = r_j+i\tau \in \gamma_j.
$$
Длина отрезка отрезка интегрирования вдоль $\gamma_j$  равна
$2hr_j^p$,   поэтому первый множитель в \eqref{Prod} оценивается
величиной
 $$
\int_0^\infty\, de(t) \left(\sup_{t\in \sigma(T)}
\sum_{j=1}^\infty\,  \frac{r_j^{2p}}
 {(t-r_j)^2} \right).
$$
Достаточно показать, что присутствующий здесь супремум ограничен
(так как $e(t)$ --- ограниченная монотонная функция). Фиксируем
число $t\in \sigma(T)$.  Найдется число $n\in \mathbb N$, такое, что
$$
t\in \Delta_n = [\alpha_{2n-1},\,\alpha_{2n}], \quad\text{а тогда}\
\,  r_{n-1} + b_1 r_{n-1}^p\leqslant t \leqslant r_n - b_1 r_n^p.
$$
Учитывая оценку \eqref{r},  получаем
\begin{align}
\nonumber \sum_{j=1}^\infty\,  \frac{r_j^{2p}} {(t-r_j)^2} &=
\frac{r_n^{2p}} {(t-r_n)^2} + \left( \sum_{j=1}^{n-1} +
\sum_{j=n+1}^\infty \right)\frac{r_j^{2p}}
{(t-r_j)^2}\\
\nonumber &\quad \leqslant  \frac{r_n^{2p}} {(r_n - b_1r^p - r_n)^2}
+ \left( \sum_{j=1}^{n-1}
  + \sum_{j=n+1}^\infty\right)\, \frac{r_j^{2p}}{(r_n-r_j)^2} \\
\label{r_n} &\quad\leqslant
   \frac 1{b_1^2} + \frac 1{c^2} \left( \sum_{j=1}^{n-1}
  + \sum_{j=n+1}^\infty\right)\, \frac 1{(n-j)^2} < C.
\end{align}
Здесь мы учли, что конечное число слагаемых в сумме при $n\leqslant
n_0$,  для которых \eqref{r}  выполняется с другой константой, на
оценку не влияет.

Остается оценить второй множитель в \eqref{Prod}. С учетом того, что
$\lambda\in \gamma_j$, а длина отрезка $\gamma_j$  равна $2hr_j^p$,
для  этого множителя находим мажоранту (с точностью до константы)
$$
\sum_{j=1}^\infty \, r_j^{2p-2} \int_0^\infty \, \frac{\ t^2 \,
dt}{(t- r_j)^2} = \int_0^\infty \, de(t) \left(\sup_{t\in \sigma(T)}
\sum_{j=1}^\infty\, \frac{r_j^{2p-2}\, t^2}{(t - r_j)^2}\right).
$$
Нужно лишь показать ограниченность супремума. Фиксируем точку $t\in
\sigma(T)$. Разобьем  сумму под знаком супремума на две части  и
оценим каждую из частей:
$$
\left(\sum_{2r_j > t} + \sum_{2r_j \leqslant t}\right)\,
\frac{r_j^{2p-2}\, t^2}{(t - r_j)^2} \leqslant  4\sum_j \,
\frac{r_j^{2p}}{(t - r_j)^2} + 4\sum_{2r_j \leqslant t} \,
r_j^{2(p-1)}.
$$
Первая сумма уже была оценена в \eqref{r_n}. Вторая сумма конечна,
но число слагаемых в ней растет при $t \to \infty$. Для оценки этой
суммы используем неравенства \eqref{1-p}. Получаем
$$
r_j^{(1-p)} \geqslant c(j-j_0) - r_0^{1-p}, \quad \text{для всех} \
\, j\geqslant j_0,
$$
где $j_0$  --- достаточно большое   фиксированное число.  Эта оценка
влечет неравенство $r_j^{(1-p)} \geqslant (c/2)j$ для всех
достаточно больших $j\geqslant j_1$,  что обеспечивает  сходимость
второго ряда. Тем самым доказана сходимость ряда \eqref{Ser1}.

{\it Шаг 5.} Теперь усложним задачу. Предположим, что выполнено
условие \eqref{m}, то есть каждый интервал $\Theta_j = (r_j-b_1
r_j^p,\, r_j+ b_1 r_j^p)$ содержит не более $m$  собственных
значений $\mu_k$  оператора $T$, пронумерованных с учетом кратности.
Тогда интегралы по кривым $\Gamma_\pm$ оцениваются так же, как и
прежде. Однако оценку интегралов по вертикальным отрезкам $\gamma_j$
нужно проводить по-другому, поскольку ограниченность
оператор-функции $M(\lambda)$ уже нельзя гарантировать. Идея этой
оценки состоит в следующем. Мы покажем, что в полосе $\Omega_r$
ширины $2\delta r^p$, определяемой \eqref{Om},  найдется полоса
$\Omega'_r$ ширины $2\delta'r^p$  при некотором $\delta' < \delta$,
не зависящем от $ r_j$,  в которой выполняется оценка
\begin{equation}\label{Fl}
\| F(\lambda)\| \leqslant C_1, \quad \text{где}\ \, F(\lambda) =
B(A-\lambda)^{-1}, \ \ \lambda \in \Omega'_j \subset \Omega_j.
\end{equation}
Более того,  мы покажем, что полосу $\Omega'_j$ можно выбрать так,
что отрезок $\Theta'_j = \Omega'_j \bigcap \mathbb R$  не содержит
собственных значений $\mu_k$  оператора $T$.  А пока предположим,
что это уже доказано.  Воспользуемся равенством
$$
 (A-\lambda)^{-1} = (T-\lambda)^{-1} (1- F(\lambda)).
$$
Тогда
\begin{multline}\label{1+C}
|(G(\lambda)x,\, x)| =  |((T-\lambda)^{-1}\, (1-F(\lambda))\,
S(\lambda)x, \, x)| \\
\leqslant(1+C_1)\|(T-\overline\lambda)^{-1} x\|\,
\|S(\lambda)x\|\leqslant (1+C_1) \|T(T-\lambda)^{-1} x\|^p\,
\|(T-\lambda)^{-1} x\|^{2-p}.
\end{multline}
Здесь мы применили оценку \eqref{Slam}  и равенство
$\|(T-\lambda)^{-1} x\| =\|(T-\overline{\lambda})^{-1}x\|,$ которое
справедливо для  самосопряженных  операторов.  Пусть  $\gamma'_j$
---  вертикальные отрезки, проходящие через точки  $r'_j$,
являющиеся серединами отрезков $\Theta'_j$ длины $2\delta' r_j^p$.
Тогда интегралы  по $\gamma'_j$  от функции \eqref{1+C} оцениваются
точно так же, как на шаге 4. Разница лишь в том, что вместо
слагаемого $1/b^2_1$  в \eqref{r_n} появится константа $
1/(\delta')^2$.

Итак, надо доказать, что существуют полосы  $\Omega'_j \subset
\Omega_j$, в которых нет собственных значений оператора $T$ и в
которых выполняется оценка \eqref{Fl}. Воспользуемся приемом Мацаева
создания {\guillemotleft}искусственной лакуны{\guillemotright}.
Рассмотрим оператор
\begin{equation}\label{K_r}
K_j: =  \sum_{\mu_k < r_j} (\mu_k - r_j^-) (\cdot,
\varphi_k)\varphi_k -  \sum_{\mu_k \geqslant r_j} (\mu_k - r_j^+)
(\cdot,
 \varphi_k)\varphi_k, \quad r_j^\pm := r_j  \pm b_1\, r_j^p,
 \end{equation}
где $\varphi_k$  --- ортонормированные собственные векторы оператора
$T$,  отвечающие собственным значениям $\mu_k$, занумерованным с
учетом кратностей, а  в сумме  участвуют только те индексы $k$, для
которых  собственные значения $\mu_k$ принадлежат отрезку
$$
\hat{\Lambda}_j =[r^+,\, r^-]  \subset \Lambda_j = (\alpha_{2j},\,
\alpha_{2j+1}).
$$
Ранг этого оператора не превышает $m$. Положим $T_j:\, = T-K_j$.
Оператор $T_j$ остается  самосопряженным,  сохраняет систему
собственных векторов $\{\varphi_k\}_1^\infty$, но меняет собственные
значения, лежащие в интервале $\Theta_j$, сдвигая их  из этого
интервала в ближайший из его концов $r_j^+$ или $r_j^-$.

Воспользуемся равенством
\begin{equation}\label{BA}
B(A-\lambda)^{-1} = B(T_j +B -\lambda)^{-1} \left(1+
L_j(\lambda)\right)^{-1}, \quad L_j(\lambda): =K_j(T_j + B
-\lambda)^{-1}.
\end{equation}
Далее,
$$
(T_j +B -\lambda)^{-1} = (T_j  -\lambda)^{-1}\left(1 +
S_j(\lambda)\right)^{-1}, \quad S_j(\lambda): = B(T_j-\lambda)^{-1}.
$$
Заметим, что при $\lambda \in\Omega_j$ (см. оценку \eqref{TT}  на
шаге 2 для оператора $T$, не имеющего собственных значений внутри
лакун $\Lambda_j$)
\begin{multline}\label{Sla}
\|S_j(\lambda)\| \leqslant b \|T(T_j -\lambda)^{-1}\|^p\, \|(T_j
-\lambda)^{-1}\|^{1-p} +M\|(T_j-\lambda)^{-1}\|  \leqslant \\
b\|TT_j^{-1}\| \|T(T_j -\lambda)^{-1}\|^p\, \|(T_j
-\lambda)^{-1}\|^{1-p} +M\|(T_j-\lambda)^{-1}\|  \leqslant\\
\frac b{b_1}\|TT_j^{-1}\|(1+o(1)) \leqslant \frac b{b_1}(1+o(1)),
\quad \lambda \in \Omega_j,
\end{multline}
так как
$$
\|TT_j^{-1}\| \leqslant \frac{r_j}{r_j^-} = \frac{r_j}{r_j +
b_1r_j^p} = 1+o(1) \quad \text{при}  \ r_j\to\infty.
$$
Но тогда
\eqref{Sl})  получаем
\begin{equation}\label{-1}
 \|(1+ S_j(\lambda))^{-1}\| \leqslant
 \frac{2(b_1 +b')}{b_1 -b')} < C, \lambda \in \Omega_j, |\lambda| \gg 1.
\end{equation}
Следовательно, в полосе $\Omega_j$ справедлива оценка
\begin{equation}\label{T_r}
\|(T_j+B-\lambda)^{-1}\| \leqslant \|(T_j - \lambda)^{-1}\|\
\|(1+S_j(\lambda))^{-1}\|   \leqslant Cr_j^{-p}, \ \ \lambda\in
\Omega_j.
\end{equation}
\tolerance=1000 Здесь мы учли простую оценку для резольвенты
самосопряженного оператора $\|(T_j(\lambda))^{-1}\| \leqslant 1/
\dist(\lambda,\, \sigma(T_j))$. Из  оценок \eqref{Sla} и  \eqref{-1}
получаем
$$
\|B(T_j+B-\lambda)^{-1}\| = \|S_j(\lambda)(1+S_j(\lambda))^{-1}\|
\leqslant C_2, \quad \lambda \in \Omega_j.
 $$

Итак, учитывая представление \eqref{BA}, для доказательства
неравенств \eqref{Fl} остается оценить оператор-функцию
\begin{equation}\label{finite}
(1+L_j^{-1}(\lambda))^{-1},  \quad    \ \, L_j(\lambda) = K_j(T_j +B
-\lambda)^{-1}.
\end{equation}
Из  определения  операторов $K_j$ следует, что $\|K_j\| \leqslant
b_1 r^p$.  Из \eqref{T_r} получаем $\|L_j(\lambda)\| \leqslant C$.
Теперь для оценки оператор-функции \eqref{finite} воспользуемся
Предложением 3, полагая  $\mu = -1, K= L_j(\lambda)$. Заметим, что в
рассматриваемом случае числитель в правой части \eqref{K}
оценивается величиной $ (1+C)^m$, так как $s$-числа оператора $K =
L_j(\lambda)$  не превышают его нормы.  Знаменатель есть скалярная
функция $ D(\lambda):= \det(1+L_j(\lambda))$ --- произведение
$\leqslant m$ чисел $1+\lambda_k(L_j(\lambda))$,  где
$\lambda_k(L_j(\lambda))$
--- собственные значения оператора$L_j(\lambda)$. Собственные
значения также оцениваются постоянной $C$ (нормой оператора
$L_j(\lambda)$),  поэтому во всей полосе $\Omega_j$ справедлива
оценка сверху $|D(\lambda)| \leqslant (1+C)^m$. В силу Предложения 2
имеем $\|(T_r+B-\lambda)^{-1}\| <C\tau^{-1},  \ \tau = \im\,
\lambda$, если $\lambda$ лежит на сторонах параболы $P_{p,h},\ h>
b'$ или вне этой параболы.   Возьмем $ h= 2Cb_1>b'.$ Тогда
$\|L_j(r_j+ihr_j^p)\|\leqslant 1/2$ и, значит,  модули всех
собственных значений оператора $L_j(r_j+ihr_j^p)$  оцениваются
сверху  $1/2$. Следовательно, имеется оценка снизу $|D(\lambda)| >
(1-1/2)^m =2^{-m}$  в точке   $\lambda = r_j^p +ihr_j^p$.

Далее воспользуемся  Предложением 4 для оценки  снизу  функции $D$.
В качестве прямоугольника $\Pi$  возьмем прямоугольник с центром в
точке $r_j$, высотой $4hr_j^p$ и шириной $2\delta r_j^p$, а в
качестве $\Pi'$ --- вдвое меньший прямоугольник с тем же центром.
Функция $D$   в прямоугольнике $\Pi$  ограничена постоянной
$M=(1+c)^m$, а для постоянной $M'$  в прямоугольнике $\Pi'$  мы
получили оценку снизу $M'\geqslant 2^{-m}$. Согласно Предложению 4
функция  $D$ оценивается снизу во всем прямоугольнике $\Pi'$ некой
постоянной вне исключительного множества кружков с общей суммой
радиусов $\varepsilon  r_j^p$. Кроме того, число нулей  функции $D$
в прямоугольнике $\Pi'$ с учетом кратности не превышает $Cm$.
Поэтому $\varepsilon$ можно выбрать столь малым, чтобы полоса
$\Pi''$ ширины  $c r^p$ в в полосе  $\Omega'$    при некотором
достаточно малом  $c>0$  не пересекала исключительное множество.
Тогда в этой полосе $\Pi''$ оператор-функция
$(1+L_j(\lambda))^{-1}$,   а потому и оператор-функция
$B(A-\lambda)^{-1}$ равномерно ограничены постоянной, не зависящей
от  $r_j \to \infty$. Этим заканчивается доказательство нужных
оценок.

{\it Шаг 6.} Нужно еще доказать полноту системы проекторов
$\{Q_j\}$.  Выберем произвольные точки $r'_n \in \mathbb R,
n\geqslant N,$  в полосах $\Pi''=\Pi''_n$, построенных на предыдущем
шаге. На прямых $\re\,\lambda = r'_n$  оператор-функции $F(\lambda)
= B(A-\lambda)^{-1}$ и $(T-\lambda)^{-1}$ равномерно ограничены.
Обозначим через $R_n$ прямоугольник, вертикальные стороны которого
проходят через точки $r'_{n-1}$ и $r'_n$, а горизонтальные стороны
таковы, что $R_n$ содержит весь спектр обоих операторов $T$ и $A$
между прямыми $\re\,\lambda = r'_{n-1}$ и  $\re\,\lambda = r'_n$.
Через $\widehat R_n$ обозначим  квадрат с центром в нуле и длиной
стороны  $2r'_n$. Из представления  \eqref{65} получаем
$$
\sum_{j=N}^n Q_jx=-\frac{1}{2\pi i}\sum_{j=N}^n \int_{\partial
R_n}\left(A-\lambda \right)^{-1}x\,\,d\lambda=\sum_{j=N}^n
P_jx+\frac{1}{2\pi i}\sum_{j=N}^n \int_{\partial R_n}G(\lambda)x\,\,
d\lambda,
$$
где $P_j$ --- ортогональные проекторы Рисса на инвариантные
подпространства оператора $T$, такие, что спектр сужения оператора
$T$ на эти инвариантные подпространства содержится в отрезках
$[r'_{j-1},\, r'_j]$. Через $P_0$  и $Q_0$  обозначим  проекторы
Рисса операторов $T$ и  $A$, отвечающие за спектры, лежащие в
полуплоскости  $\re\,\lambda < r'_N$. Далее для упрощения записи
считаем $N=1.$ Положим
$$
 \widehat P_n = \sum_{j=0}^n P_j, \qquad \widehat Q_n = \sum_{j=0}^n Q_j,
 \qquad H_n = \hat{P_n}(H).
 $$
Тогда
\begin{equation}\label{Gn}
\widehat Q_nx = \widehat P_n x + \widehat G_n x, \qquad \text{где}\
\widehat G_n = \frac{1}{2\pi i}\int_{\partial \widehat R_n}\,
G(\lambda)\, d\lambda.
\end{equation}
Покажем, что
\begin{equation}\label{Gla}
 \widehat G_nx: =\int_{\partial \widehat R_n}\, G(\lambda)x\,
 d\lambda \to 0 \quad  \text{при}\, n\to \infty \  \,\forall\, x\in H.
\end{equation}
Тогда
$$
 \sum_{j=0}^n Q_jx \to x \ \, \text{при}\ n\to\infty, \ \,
 \text{поскольку, очевидно,}\ \sum_{j=0}^n P_jx\to x\ \text{при}\  n\to\infty.
$$
Первое соотношение влечет полноту системы проекторов $\{Q_j\}$,
поэтому достаточно доказать  соотношение \eqref{Gla}.

Пусть  $e(t) = (E(t)x,x)$, где $E(t)$ --- спектральная функция
оператора $T$. Для всех $x\in H_n$ имеем
$$
 \left\|\left(T-\lambda \right)^{-1}x_n\right\|^2=
 \int\limits_{-r'_n}^{r'_n}\frac{de(t)}{(\xi -t)^2+\tau^2}, \quad
 \lambda=\xi+i\tau,\quad \xi,\tau\in\mathbb{R},
$$
откуда получаем (при фиксированном $n$)
\begin{equation}\label{small}
\left\|\left(T-\lambda \right)^{-1}x_n\right\|^2=
O\left(|\lambda|^{-1}\right)\ \ \ \text{при}\ \,  |\lambda|\to
\infty.
\end{equation}
Не ограничивая общности, будем считать оператор $T$ обратимым (иначе
можно провести сдвиг спектрального параметра). Так как $H_n$
инвариантно относительно $T$, то оператор $T$  изоморфно отображает
$H_n$ на  $H_n$. Поэтому для всех $x_n\in H_n$
\begin{equation}\label{infty}
 G(\lambda)x_n = (T-\lambda)^{-1} M(\lambda)BT^{-1}(T-\lambda)^{-1}y_n,
 \quad \text{где } \  y_n = Tx_n.
\end{equation}
Так как оператор-функция $M(\lambda)BT^{-1}$ ограничена на сторонах
квадрата $\widehat{R}_n$,  то в силу \eqref{small} и \eqref{infty}
получаем (при фиксированном $n$)
\begin{equation}\label{j}
\int_{\partial \widehat R_j}\,  G(\lambda) x_n\, d\lambda \to
0,\quad \text{если} \  j\to \infty,
\end{equation}
так как на горизонтальных и левой сторонах квадрата $\widehat R_j$
имеем $\|G(\lambda)x_n\| = O(|\lambda|^{-2})$, а на правой стороне
$\|G(\lambda)x_n\| = O(|\lambda|^{-1-p})$ (при $p=0$  соотношение
\eqref{j} также выполняется, если учесть оценку
$\|(T-\lambda)^{-1}\| \leqslant |\im\,\lambda|^{-1}$).

Таким образом соотношение \eqref{Gla}  имеет место на плотном
множестве в $H$ (поскольку  множество векторов $x_n\in H_n$ при всех
$n$ плотно в $H$). Поэтому для доказательства соотношения
\eqref{Gla} достаточно показать, что  нормы $\|\widehat G_n\|$
равномерно ограничены. Из сходимости ряда \eqref{67}  следует
равномерная ограниченность $|(\widehat Q_n x,\, x)|$ при всех  $x\in
H$. Квадратичная форма $|(\widehat Q_n x,\, x)|$  в гильбертовом
пространстве определяет билинейную форму, поэтому $|(\widehat Q_n
x,\, y)| <\infty$ для всех $x,y\in H$. Дважды применяя теорему
Банаха-Штейнгауза, получим сначала равномерную ограниченность норм
$\|\widehat Q_n x\|$, а затем равномерную ограниченность $\|\widehat
Q_n\|$. Тогда в силу  равенства \eqref{Gn} нормы $\|\widehat G_n\|$
также равномерно ограничены. Этим заканчивается доказательство
полноты проекторов $\{Q_n\}$.

{\it Шаг 7.} Мы рассмотрели случай $T=T^* >0$.  Это условие
эквивалентно тому, что $T$ --- нормальный оператор со спектром на
одном луче. Случай, когда  оператор $T$ не является
полуограниченным, не вносит трудностей. В этом случае, согласно
Предложению 2, спектр возмущенного оператора $T+B$  лежит внутри
объединения двойной параболы $P_{p, h}$ и круга достаточно большого
радиуса. Интегралы по границам этих двойных парабол оцениваются без
изменений. Конструкция отрезков $\gamma'_j$, соединяющих ветви левой
параболы, и оценки интегралов по этим отрезкам проводятся
аналогично. Все остальные рассуждения остаются прежними, только
индексацию нужно проводить не по натуральным, а по целым числам.

\newpage

\small
\begin{flushleft}
Alexander Konstantinovich Motovilov \\
Bogoliubov Laboratory of Theoretical Physics\\ 
Joint Institute for Nuclear Research \\ 
Joliot-Curie 6\\ 
141980 Dubna, Moscow Region, Russia \\ 
E-mail: motovilv@theor.jinr.ru
\smallskip

and
\smallskip

Faculty of Natural and Engineering Sciences\\
Dubna State University\\
Universitetskaya 19\\
141980 Dubna, Moscow Region, Russia

\end{flushleft}
\bigskip

\begin{flushleft}
Andrei Andreevich Shkalikov\\ 
Faculty of Mathematics and Mechanics\\ 
Lomonosov Moscow  State University \\ 
Leninskiye Gory 1\\ 
119991 Moscow GSP-1, Russia\\ 
E-mail: shkalikov@mi.ras.ru
\end{flushleft}

\end{document}